\documentclass[11pt,letterpaper]{amsart}

\usepackage[T1]{fontenc}
\usepackage{lmodern}
\usepackage{hyperref}
\usepackage{etex}
\usepackage[shortlabels]{enumitem}
\usepackage{amsmath}
\usepackage{amsxtra}
\usepackage{amscd}
\usepackage{amsthm}
\usepackage{amsfonts}
\usepackage{amssymb}
\usepackage{eucal}
\usepackage[all]{xy}
\usepackage{graphicx}
\usepackage{tikz-cd}
\usepackage{mathrsfs}
\usepackage{subfiles}
\usepackage{mathpazo}
\usepackage[colorinlistoftodos, textsize=tiny]{todonotes}
\usepackage{morefloats}
\usepackage{pdfpages}
\usepackage{thm-restate}
\usepackage[utf8]{inputenc}
\usepackage{epigraph}
\usepackage{csquotes}
\usepackage[margin=1in]{geometry}
\usepackage{adjustbox}
\newcommand\reallywidehat[1]{%
\savestack{\tmpbox}{\stretchto{%
  \scaleto{%
    \scalerel*[\widthof{\ensuremath{#1}}]{\kern-.6pt\bigwedge\kern-.6pt}%
    {\rule[-\textheight/2]{1ex}{\textheight}}
  }{\textheight}%
}{0.5ex}}%
\stackon[1pt]{#1}{\tmpbox}%
}
\graphicspath{ {images/} }

\RequirePackage{color}
\definecolor{myred}{rgb}{0.75,0,0}
\definecolor{mygreen}{rgb}{0,0.5,0}
\definecolor{myblue}{rgb}{0,0,0.65}

\usepackage{color}

\usepackage{hyperref}
\hypersetup{citecolor=blue}
\usepackage{tikz}
\usetikzlibrary{matrix,arrows,decorations.pathmorphing}

\theoremstyle{plain}
\newtheorem{theorem}[subsection]{Theorem}

\newtheorem{proposition}[subsection]{Proposition}
\newtheorem{lemma}[subsection]{Lemma}
\newtheorem{corollary}[subsection]{Corollary}

\theoremstyle{definition}
\newtheorem{definition}[subsection]{Definition}
\newtheorem{remark}[subsection]{Remark}

\theoremstyle{remark}

\numberwithin{equation}{subsection}
\newcommand\nc{\newcommand}
\nc\on{\operatorname}
\nc\renc{\renewcommand}

\newcommand\bc{\mathbb C}

\newcommand\bl{\mathbb L}

\newcommand\bp{\mathbb P}
\newcommand\bq{\mathbb Q}

\newcommand\bv{\mathbb V}

\newcommand\bz{\mathbb Z}

\newcommand*{\shom}{\mathscr{H}\kern -.5pt om}
\newcommand*{\stor}{\mathscr{T}\kern -.5pt or}
\newcommand*{\sext}{\mathscr{E}\kern -.5pt xt}

\makeatletter
\providecommand\@dotsep{5}
\renewcommand{\listoftodos}[1][\@todonotes@todolistname]{%
\@starttoc{tdo}{#1}}
\makeatother

\makeatletter
\newcommand{\customlabel}[2]{\protected@write \@auxout {}{\string \newlabel {#1}{{#2}{\thepage}{#2}{#1}{}} }\hypertarget{#1}{#2}}

\DeclareMathOperator\id{id}

\DeclareMathOperator\rank{rank}
\DeclareMathOperator\codim{codim}

\DeclareMathOperator\gal{Gal}

\DeclareMathOperator\gl{GL}

\renewcommand\u{\mathrm{U}}

\DeclareMathOperator\pardeg{par-deg}

\DeclareFontFamily{U}{wncy}{}
\DeclareFontShape{U}{wncy}{m}{n}{<->wncyr10}{}
\DeclareSymbolFont{mcy}{U}{wncy}{m}{n}
\DeclareMathSymbol{\Sha}{\mathord}{mcy}{"58}

\def\listtodoname{List of Todos}
\def\listoftodos{\@starttoc{tdo}\listtodoname}

\title{Hypergeometric local systems and parabolic bundles}

\author{Charlie Wu}

\usepackage{microtype}

\begin{document}

\begin{abstract}
Beukers and Heckman gave necessary and sufficient conditions for a hypergeometric function $_n F_{n-1}$ to be algebraic. We give a new proof of this theorem by passing through the Mehta-Seshadri correspondence. In particular, we explicitly write down the parabolic bundle corresponding to a unitary hypergeometric local system.
\end{abstract}
\maketitle

\addtocontents{toc}{\protect\setcounter{tocdepth}{1}}
\vspace{-.5cm}

\setcounter{tocdepth}{1}
\tableofcontents

\vspace{-.5cm}
\section{Introduction}\label{sec:section 1}
For a number $q$, let $(q)_k = q(q + 1) \cdots (q + k - 1)$ be its Pochhammer symbol. For a tuple of numbers $(\alpha, \beta) = (\alpha_1, \dots, \alpha_n, \beta_1, \dots, \beta_{n-1})$, the hypergeometric function $_n F_{n-1}(\alpha, \beta, z)$ is defined to be
\begin{align*}
	_n F_{n-1}(\alpha_1, \dots, \alpha_n, \beta_1, \dots, \beta_{n-1} , z) &= \sum_{k = 0}^\infty \frac{(\alpha_1)_k \cdots (\alpha_n)_k}{(\beta_1)_k \cdots (\beta_{n-1})_k \cdot k!} z^k.
\end{align*}
The hypergeometric function $_n F_{n-1}$ satisfies a differential equation $P(\alpha, \beta)s = 0$ on $\bp^1$ where
\begin{align*}
P(\alpha_1, \dots, \alpha_n, \beta_1,\dots, \beta_{n-1}) = D(D + \beta_1 - 1) \cdots (D + \beta_{n - 1} - 1) - z (D + \alpha_1) \cdots (D + \alpha_n)
\end{align*}
and $D = z \frac{d}{dz}$. This differential equation has regular singularities at $0$, $1$, and $\infty$. We call the corresponding local system on $\bp^1 \setminus \{0,1,\infty\}$ under the Riemann-Hilbert correspondence a \emph{hypergeometric local system}. This local system has local sections $s$ given by solutions to the differential equation $P(\alpha, \beta)s = 0$.

Let $\bv$ be a local system on a smooth variety $S$. We say that $\bv$ is \emph{of geometric origin} if there is a Zariski-dense open $U \subseteq S$ and a smooth proper morphism $f: X \to U$ so that $\bv|_U$ is a subquotient of $R^i f_* \bc$ for some $i \geq 0$. When the parameters $\alpha_i$ and $\beta_i$ are rational numbers, the hypergeometric local systems are all of geometric origin. To see this, one expresses $_n F_{n-1}$ as a period integral \cite[Equation (4.1.3)]{Slater_generalized}

\begin{align*}
	_n F_{n-1}(\alpha, \beta, z) &= \prod_{j = 1}^{n-1}\frac{\Gamma(\beta_{j})}{\Gamma(\alpha_{j+1})\Gamma(\beta_{j} - \alpha_{j+1})}\int_0^1 \cdots \int_0^1 \prod_{j = 1}^{n - 1} t_j^{\alpha_{j + 1} - 1}(1 - t_j)^{\beta_j - \alpha_{j + 1}- 1} (1 - z t_1 \cdots t_{n-1})^{- \alpha_1} \,dt_1 \dots \,dt_{n-1}.
\end{align*}

Monodromy groups of hypergeometric local systems have been studied extensively, especially in the work of Beukers and Heckman. Most significantly, they determine when the monodromy group is finite \cite[Theorem 4.8]{Beukers_Heckman} in terms of the parameters $\alpha_i$ and $\beta_i$. Equivalently, they determine when $_n F_{n-1}$ is an \emph{algebraic} function.  Their condition states that $_n F_{n-1}$ is algebraic if and only if all of the $\alpha_i$ and $\beta_i$ are rational numbers and all of the Galois conjugates of the sets $\{e^{2\pi i \alpha_1}, \dots, e^{2\pi i \alpha_n} \}$ and $\{ e^{2\pi i \beta_1}, \dots, e^{2\pi i \beta_{n-1}}, 1 \}$ \emph{interlace} in the sense of \autoref{def:interlacing}.

 Another proof of the theorem of Beukers and Heckman is given by Fedorov \cite[Corollary 2.9]{Fedorov}. Fedorov studies hypergeometric local systems using Katz's middle convolution operation \cite{Katz_rigid_local_systems}. Beukers and Heckman, and later Fedorov, show finiteness of the monodromy of $\mathcal{H}_n$ by showing first that it is unitary. Beukers and Heckman show this by constructing a positive-definite Hermitian form that the monodromy preserves, and Fedorov shows this by computing the Hodge numbers of the hypergeometric local system and showing that the criterion given by Beukers and Heckman gives a positive-definite flat Hermitian metric on the polarizable variation of Hodge structure associated to the local system via the Riemann-Hilbert correspondence. 

We give an alternative proof of Beukers' and Heckman's theorem by passing through the Mehta-Seshadri correspondence \cite[Theorem 4.1]{Mehta-Seshadri}, which gives a correspondence between unitary local systems on curves and stable parabolic bundles. Our main theorem gives an explicit description of the parabolic bundle associated to a unitary hypergeometric local system.
\begin{theorem}
	 There exists a stable parabolic bundle $E_\star$ corresponding to a hypergeometric local system  via the Mehta-Seshadri correspondence (equivalently, the hypergeometric local system is unitary) if and only if the sets $\{e^{2\pi i \alpha_1}, \dots, e^{2\pi i \alpha_n}\}$ and $\{e^{2\pi i \beta_1}, \dots, e^{2\pi i \beta_{n-1}}, 1\}$ interlace in the sense of \autoref{def:interlacing}. When the parameters do interlace, we give an explicit description of $E_\star$.
\begin{proof}
	See \autoref{thm:underlying bundle} and \autoref{thm:unitary iff interlacing}.
\end{proof}
 \end{theorem}
We will see that the stability conditions for parabolic bundles naturally give the conditions found in the original proof of Beukers and Heckman.

\subsection{Structure of the paper}
In \autoref{sec:section 2}, we give background on parabolic bundles and the Mehta-Seshadri correspondence. In \autoref{sec:section 3}, we recall some basic results on hypergeometric local systems. We give the alternative proof of the theorem of Beukers and Heckman in \autoref{sec:section 4}.

\subsection{Acknowledgements}
The author would like to thank Daniel Litt for many helpful discussions. The author would also like to thank Daniel Litt and Matthew Bolan for running a seminar on hypergeometric functions in which the author first learned about hypergeometric local systems and the results of Beukers and Heckman.

\section{Background on parabolic bundles and the Mehta-Seshadri correspondence}\label{sec:section 2}

We review the Mehta-Seshadri correspondence \cite[Theorem 4.1]{Mehta-Seshadri} in this section. This correspondence is crucial to our study of hypergeometric local systems with finite monodromy.

Let $X$ be a smooth proper curve over $\bc$ and $D = x_1 + \dots + x_n$ a reduced effective divisor on $X$.

\begin{definition}\label{def: parabolic structure}
A \emph{parabolic vector bundle} on $(X,D)$ is a vector bundle $E$ on $X$ together with the following data:
\begin{enumerate}
	\item a strictly decreasing filtration of subspaces of the fiber $E_{x_j} = E^1_j \supset E^2_j \supset \dots \supset E^{n_j}_j \supset E^{n_j + 1}_j = 0$ for each $x_j \in D$,
	\item a sequence of real numbers $0 \leq \alpha_j^1 < \alpha_j^2 < \dots < \alpha_j^{n_j} < 1$ for each $x_j \in D$.
\end{enumerate}
We write $E_\star$ to mean the data $(E, \{E_j^i\}, \{\alpha_j^i\})$. We call the $\{\alpha_j^i\}$ the \emph{parabolic weights} of $E_\star$.
\end{definition}

\begin{definition}\label{def:parabolic degree and slope}
Let $E_\star$ be a parabolic bundle on $(X,D)$. We define the \emph{parabolic degree} to be
\begin{align*}
	\pardeg E_\star &:= \deg E + \sum_{j = 1}^n \sum_{i = 1}^{n_j} \alpha_j^i \cdot \dim(E_j^i/E_j^{i + 1}).
\end{align*}
We define the \emph{parabolic slope} to be $\mu_\star(E_\star) := \pardeg E_\star/ \rank(E_\star)$.
\end{definition}

\begin{definition}\label{def:parabolic structure for subquotients}
Let $E_\star$ be a parabolic bundle and $F$ a subbundle of $E$. Both $F$ and $E/F$ have induced parabolic structures. For $F$, the filtration of the fiber $F_{x_j}$ is obtained by removing redundancies from
\begin{align*}
	F_{x_j} = E^1_j \cap F_{x_j} \supseteq \dots \supseteq E^{n_j}_j \cap F_{x_j} \supseteq E^{n_j +1 }_j \cap F_{x_j} = 0.
\end{align*}
The parabolic weight associated to $F_j^i$ is given by $\max_{k, 1 \leq k \leq n_j}\{\alpha_j^k :F_j^i = E_{j}^k \cap F_{x_j}\}$.
Similarly, the filtration defining the parabolic structure on $E/F$ is obtained by removing redundancies from
\begin{align*}
	(E/F)_{x_j} = (E_j^1 + F_{x_j})/F_{x_j} \supseteq (E_j^2 + F_{x_j})/F_{x_j} \dots \supseteq (E_j^{n_j} + F_{x_j})/F_{x_j} \supseteq (E_j^{n_j + 1} + F_{x_j})/F_{x_j} = 0
\end{align*}
for each $x_j \in D$. The parabolic weight associated to $(E/F)_j^i$ is given by $\max_{k, 1 \leq k \leq n_j}\{\alpha_j^k :(E/F)_j^i = (E_j^k + F_{x_j})/F_{x_j}\}$.
\end{definition}

\begin{definition}\label{def:stability}
We say that a parabolic bundle $E_\star$ is \emph{semistable} if for all subbundles $F_\star$ of $E_\star$, the slopes satisfy the inequality $\mu_\star(F_\star) \leq \mu_\star (E_\star)$. We say that $E_\star$ is \emph{stable} if for all nonzero proper subbundles $F_\star$ of $E_\star$, $\mu_\star (F_\star) < \mu_\star (E_\star)$. We say that $E_\star$ is polystable if it is isomorphic to a direct sum of stable parabolic bundles of the same parabolic slope.
\end{definition}

We can now state the Mehta-Seshadri correspondence.
\begin{theorem}\label{thm:Mehta-Seshadri}
Let $C_1,\dots, C_n \subseteq \u(r)$ be conjugacy classes with $C_j$ represented by the diagonal matrix with eigenvalues $e^{2\pi i \alpha_j^1}, \dots, e^{2\pi i \alpha_j^{n_j}}$ with multiplicities $m_j^1,\dots, m_j^{n_j}$. Then, there is a bijection between (isomorphism classes of) unitary local systems on $X \setminus D$ with local monodromy $C_1,\dots, C_n$ and (isomorphism classes of) polystable parabolic bundles $E_\star$ of parabolic degree zero on $(X,D)$ with weights $\alpha_j^i$ at $x_j$ and flags $E_j^\bullet$ of $E_{x_j}$ such that $m_j^{s} = \dim(E_j^s/ E_j^{s + 1})$. Under this bijection, the irreducible unitary local systems correspond to the stable parabolic bundles of parabolic degree zero.
\begin{proof}
	See \cite[Theorem 4.1]{Mehta-Seshadri}.
\end{proof}
\end{theorem}

\section{Basic results on hypergeometric local systems}\label{sec:section 3}

Let $(\alpha,\beta) = (\alpha_1,\dots, \alpha_n,\beta_1,\dots,\beta_n)$ and let $\mathcal{H}_n(\alpha, \beta)$ be the local system on $\bp^1 \setminus \{0,1,\infty\}$ associated to the hypergeometric differential equation
\begin{align*}
	(D + \beta_1 - 1) \cdots (D + \beta_n - 1) - z (D + \alpha_1) \cdots (D + \alpha_n) = 0.
\end{align*}
We will often write $\mathcal{H}_n$ for $\mathcal{H}_n(\alpha,\beta)$ when the context is clear. We recall some basic facts about $\mathcal{H}_n$ in this section.

\begin{remark}
This differential equation is a generalization of the differential equation satisfied by the hypergeometric function $_n F_{n-1}$. Setting $\beta_1 = 1$ recovers the original differential equation.
\end{remark}

\begin{proposition}\label{prop:hypergeometric irreducible}
	The local system $\mathcal{H}_n$ is irreducible if and only if $\alpha_j - \beta_k \not \in \bz$ for all $j$ and $k$.
	\begin{proof}
		See \cite[Proposition 2.5 and Proposition 2.7]{Beukers_Heckman}.
	\end{proof}
\end{proposition}

\begin{proposition}\label{prop:defined over ring of integers}
	Suppose that all of the parameters $\alpha_i$ and $\beta_i$ are rational numbers. Let $K$ be the number field $\bq(e^{2\pi i \alpha_1}, \dots, e^{2\pi i \alpha_n}, e^{2\pi i \beta_1}, \dots, e^{2\pi i \beta_n})$. Then, $\mathcal{H}_n$ is defined over the ring of integers $\mathcal{O}_K$.
	\begin{proof}
		See \cite[Theorem 1.1]{Levelt} or \cite[Theorem 3.5]{Beukers_Heckman}.
	\end{proof}
\end{proposition}

We describe the local monodromy of $\mathcal{H}_n$, and we show that $\mathcal{H}_n$ is determined by its local monodromy.
\begin{definition}\label{def:pseudoreflection}
	Let $A$ be an invertible matrix. We say that $A$ is a \emph{pseudoreflection} if $\rank (A - I) = 1$. We call $c = \det(A)$ the \emph{special eigenvalue} of $A$.
\end{definition}

\begin{proposition}\label{prop:local monodromy}
The local monodromies of $\mathcal{H}_n$ at $0$, $1$, and $\infty$ are (conjugate to)
\begin{enumerate}
	\item a diagonal matrix with eigenvalues $e^{2\pi i \alpha_1},\dots, e^{2\pi i \alpha_n}$ around $0$,
	\item a diagonal matrix with eigenvalues $e^{- 2\pi i \beta_1}, \cdots, e^{- 2\pi i \beta_n}$ around $\infty$,
	\item and a pseudoreflection around $1$ with special eigenvalue $c = e^{2\pi i \gamma}$ where $\gamma = \sum_{i = 1}^n (\beta_i - \alpha_i)$.
\end{enumerate}
\begin{proof}
	See \cite[Proposition 3.2]{Beukers_Heckman}.
\end{proof}
\end{proposition}

\begin{definition}\label{def:rigid local system}
	Let $\mathbb{V}$ be a rank $n$ local system on $\bp^1\setminus \{x_1,\dots, x_m\}$.	
	For each $i = 1,\dots, m$, let $C_i$ be the conjugacy class in $\gl_n(\bc)$ containing the local monodromy of $\bv$ around the puncture $x_i$. We say that the local system $\mathbb{V}$ is \emph{rigid} if any other local system with local monodromy data lying in $C_i$ around each puncture $x_i$ is isomorphic to $\mathbb{V}$.
\end{definition}

\begin{proposition}\label{prop:hypergeometric is rigid}
The local system $\mathcal{H}_n$ is rigid.
\begin{proof}
	To show this, we use Katz's criterion for rigidity \cite[Theorem 1.1.2]{Katz_rigid_local_systems}. Let $\mathbb{V}$ be an irreducible rank $n$ local system on $\bp^1 \setminus \{x_1,\dots, x_m\}$ where the homotopy class of the loop $\gamma_i$ around $x_i$ is sent to the matrix $A_i \in \gl_n(\bc)$. Katz's criterion states that $\mathbb{V}$ is a rigid local system if and only if
	\begin{align*}
		2 = (2 - m)n^2 + \sum_{i = 1}^m \dim Z(A_i)
	\end{align*}
	where $Z(A_i)$ is the centralizer of $A_i$ in $\gl_n(\bc)$. In the case where $\mathbb{V}$ is a hypergeometric local system $\mathcal{H}_n$, $m = 3$ and around $0$ and $\infty$ our matrices are conjugate to diagonal matrices with $n$ distinct eigenvalues. Around $1$, our matrix is conjugate to a pseudoreflection. For a diagonal matrix $A$ with distinct eigenvalues, $\dim Z(A) = n$. For a pseudoreflection $R$, $\dim Z(R) = (n - 1)^2 + 1$. By substituting the above into the expression given in Katz's criterion, we obtain
	\begin{align*}
		(2 - 3)n^2 + \sum_{i = 1}^3 \dim Z(A_i) = -n^2 + 2n + (n - 1)^2 + 1 = 2,
	\end{align*}
	giving us rigidity of the local system $\mathcal{H}_n$.
\end{proof}
\end{proposition}

\section{Finite monodromy}\label{sec:section 4}

In this section, we explicitly write down the parabolic bundle associated to a \emph{unitary} hypergeometric local system. Using this explicit description, we give a new proof of a classical result of Beukers and Heckman characterizing when the hypergeometric local system has finite monodromy \cite[Theorem 4.8]{Beukers_Heckman}. The most important property of $\mathcal{H}_n$ that we use is that it is rigid (\autoref{def:rigid local system}). This will allow us to conclude that \emph{any} polystable parabolic bundle with the desired parabolic weights corresponds to $\mathcal{H}_n$ under the Mehta-Seshadri correspondence (\autoref{thm:Mehta-Seshadri}).

We first observe some easy necessary conditions for $\mathcal{H}_n$ to have unitary monodromy. Unitary matrices are diagonalizable with eigenvalues lying on the unit circle, so the parameters $(\alpha, \beta) = (\alpha_1, \dots, \alpha_n, \beta_1, \dots, \beta_n)$ must be tuples of real numbers. All of the $\alpha_i$ (and all of the $\beta_i$) must be distinct, since otherwise the local monodromy matrices at $0$ and $\infty$ will not be diagonalizable. The local system is the same after adding an integer to each $\alpha_i$ and $\beta_i$ (due to rigidity of $\mathcal{H}_n$), so we may assume that all of the $\alpha_i$ and $\beta_i$ lie in $[0,1)$. We make the assumption throughout this section that $\mathcal{H}_n$ is irreducible, which is the same as saying the $\alpha_i$ and $\beta_j$ are never equal (\autoref{prop:hypergeometric irreducible}). This is the same as saying the parabolic bundle corresponding to $\mathcal{H}_n$ is stable.

\subsection{Definition of the parabolic bundle $E_\star$}
We let $E_\star := E_\star(\alpha, \beta)$ be the parabolic bundle on $(\bp^1, 0 + 1 + \infty)$ given by the following data:
\begin{enumerate}
	\item The underlying bundle $E(\alpha, \beta)$ is isomorphic to $\mathscr{O}(-1)^{\oplus n}$.
	\item The flags $T^\bullet$ at $0$ and $S^\bullet$ at $\infty$ are generic \emph{complete} flags of the fibers $E_0$ and $E_\infty$, respectively. The flag $L^\bullet$ at $1$ is given by $E_1 =L^0 \supseteq L^1 \supseteq L^2 = 0$ where $L^1$ is a generic line.
	\item The parabolic weights at $0$ are $0 \leq \alpha_{i_1}< \dots < \alpha_{i_n} < 1$ where we reorder the $\alpha$'s to be in increasing order. Similarly, the parabolic weights at $\infty$ are $0 \leq 1 - \beta_{j_1} < \dots < 1 - \beta_{j_n} < 1$ where we reorder the $\beta$'s so the $1 - \beta_{j_k}$'s are in increasing order, and the parabolic weights at $1$ are $0 \leq \sum (\beta_i - \alpha_i) - \lfloor \sum (\beta_i - \alpha_i) \rfloor < 1$.
\end{enumerate}
We write $E_\star(1 - \alpha, 1 - \beta)$ to mean $E_\star(1 - \alpha_1,\dots, 1 - \alpha_n, 1 - \beta_1,\dots, 1 - \beta_n)$. If $E_\star(\alpha, \beta)$ is stable, of parabolic degree $0$, and $0 \leq \sum(\beta_i - \alpha_i)$, then $\sum(\beta_i - \alpha_i) < 1$. This is because $\pardeg E_\star(\alpha, \beta) = \lfloor \sum (\beta_i - \alpha_i) \rfloor = 0$ since $\deg E(\alpha,\beta) = -n$. For similar reasons, if $E_\star(1 - \alpha, 1 - \beta)$ is stable, of parabolic degree $0$, and $\sum(\beta_i - \alpha_i) < 0$ then $\sum(\alpha_i - \beta_i) < 1$.

\begin{remark}\label{rem:bundle well defined}
The isomorphism class of the parabolic bundle does not depend on the choice of generic flags at $0$, $1$, and $\infty$. An isomorphism of parabolic bundles consists of an isomorphism of the underlying bundles $f: E^1 \to E^2$ (which is given by any matrix $A \in \gl_n(\bc)$) so that the isomorphism sends the parabolic flags of $E^1_\star$ over each point in the parabolic divisor to the parabolic flags of $E^2_\star$. For the existence of such an $A$, see \cite[Theorem 2.2 Case $(S_{q,r})$]{Magyar-Weyman-Zelevinsky} and observe that the finite orbit condition implies there is a unique top-dimensional orbit of the diagonal action of $\gl_n(\bc)$ on $F \times F \times \bp^{n-1}$ where $F$ is the complete flag variety on $\bc^n$.
\end{remark}

\subsection{Proofs of \autoref{thm:underlying bundle}, \autoref{thm:unitary iff interlacing}, and finite monodromy} 

We write $V = H^0(\bp^1, \shom(\mathscr{O}(-1), E))$. As before, we let the parabolic flag at $0$ be given by $E_0 = T^0 \supset T^1 \supseteq \dots \supset T^{n - 1} \supset T^n = 0$. Let the parabolic flag at $\infty$ be given by $E_\infty = S^0 \supset S^1 \supset \dots \supset S^{n - 1} \supset S^n = 0$. The flag at $1$ is given by just a single line $E_1 = L^0 \supset L^1 \supseteq L^2 = 0$. Since our flags $T^\bullet$, $S^\bullet$, and $L^\bullet$ are generic, we have that $V = T^j \oplus S^{n - j} = T^j \oplus S^{n - j + 1} \oplus L^1$ for all $j  \in \{1,\dots, n\}$. 

Let $g: G^k_\star \hookrightarrow E_\star$ be an arbitrary rank $k$ subbundle. For a point $p \in \bp^1$, we write $G^k_p$ for the fiber of $G^k$ at $p$. We first prove a technical lemma about $G^k$. This lemma will be used to prove that $\pardeg G^k_\star < 0$ in the proof of \autoref{thm:unitary iff interlacing}.

\begin{lemma}\label{lem:inequalities}
Let $g: G^k_\star \hookrightarrow E_\star$ be a rank $k$ parabolic subbundle as above. Then,
\begin{align*}
\deg G^k + \dim (G^k_0 \cap T^j) + \dim(G^k_\infty \cap S^{n - j}) &\leq 0,\\
\deg G^k + \dim (G^k_0 \cap T^j) + \dim(G^k_\infty \cap S^{n - j + 1}) + \dim (G^k_1 \cap L^1) &\leq 0.
\end{align*}
\begin{proof}
We prove the second inequality. The first inequality is completely analogous.

Since all of the fibers of $E$ are canonically isomorphic to $V$, the splitting of $V$ above gives us a splitting of $E$ into $\mathscr{O}(-1) \otimes_\bc T^j \oplus \mathscr{O}(-1) \otimes_\bc S^{n - j + 1} \oplus \mathscr{O}(-1) \otimes_\bc L^1$. Using this splitting, we define the map $\varphi^j: E \to \mathscr{O}^{\oplus n}$ given by the $n \times n$ matrix
\begin{align*}
\begin{pmatrix}
x I_{n - j} & 0 & 0\\
0 & y I_{j - 1} & 0\\
0 & 0 & x - y
\end{pmatrix}.
\end{align*}
Observe that at $0$, the induced map $\varphi^j_0$ has kernel $T^j$. At $\infty$, $\varphi^j_\infty$ has kernel $S^{n - j + 1}$. And at $1$, $\varphi^j_1$ has kernel $L^1$. And away from $0$, $1$, and $\infty$ this map is invertible.

We then restrict this map to $G^k$. Observe that $G^k_0 \cap T^j$ is the kernel of $\varphi^j|_{G^k}$ at $0$, $G^k_1 \cap L^1$ is the kernel of $\varphi^j|_{G^k}$ at $1$, and $G^k_\infty \cap S^{n - j + 1}$ is the kernel of $\varphi^j|_{G^k}$ at $\infty$. 

Let $p \in \bp^1$ be a closed point which is not $0$, $1$, or $\infty$. Then at $p$, the inclusion $g: G^k_\star \hookrightarrow E_\star$ gives us an inclusion of fibers $g_p: G^k_p \hookrightarrow V$. Let $g_p(G^k_p)$ be the image of this map inside $V$. Then, $\varphi_p^j(g_p(G^k_p))$ is a $k$-dimensional subspace of $V$. We let $\pi: V \to \varphi_p^j(g_p(G^k_p))$ be a projection of $V$ onto this subspace so that restricting $\pi$ to $\varphi_p^j(g_p(G^k_p))$ yields an isomorphism.

Now, consider the composition 
\begin{align*}
G^k \xrightarrow{g} V \otimes \mathscr{O}(-1) \xrightarrow{\varphi^j} V \otimes_\bc \mathscr{O} \xrightarrow{\pi \otimes \id} (\varphi^j \circ g)_p(G^k_p) \otimes_\bc \mathscr{O}.
\end{align*}
We write $h : G^k \to (\varphi^j \circ g)_p(G^k_p) \otimes_\bc \mathscr{O}$ to be the composition of the maps above. Note that $(\varphi^j \circ g)_p(G^k_p) \otimes_\bc \mathscr{O}$ is isomorphic to $\mathscr{O}^{\oplus k}$. Then, we obtain a nonzero global section $\det h \in H^0(\bp^1, \det( G^k)^\vee)$ because at $p$, $h$ is an isomorphism so $\det h$ is nonzero at $p$. Furthermore, at $0$, $\infty$, and $1$, $\det h$ vanishes to order at least $\dim(G^k_0 \cap T^j)$, $\dim(G^k_\infty \cap S^{n - j + 1})$, and $\dim(G^k_1 \cap L^1)$, respectively. This is because
\begin{align*}
G^k_0 \cap T^j \subseteq &\ker h_0\\
G^k_\infty \cap S^{n - j + 1} \subseteq & \ker h_\infty\\
G^k_1 \cap L^1 \subseteq &\ker h_1.
\end{align*}

Let $\mathrm{div} (\det h)$ be the divisor on $\bp^1$ cut out by $\det h$. Then,
\begin{align*}
\deg (\det (G^k)^\vee) = - \deg G^k = \deg (\mathrm{div} (\det h)) \geq \dim(G^k_0 \cap T^j) +\dim(G^k_\infty \cap S^{n - j + 1}) +  \dim(G^k_1 \cap L^1).
\end{align*}
In other words,
\begin{align*}
\deg G^k + \dim(G^k_0 \cap T^j) +\dim(G^k_\infty \cap S^{n - j + 1}) +  \dim(G^k_1 \cap L^1) \leq 0,
\end{align*}
as desired.
\end{proof}
\end{lemma}

\begin{theorem}\label{thm:underlying bundle}
If $\sum(\beta_i - \alpha_i) \geq 0$, then $\mathcal{H}_n(\alpha, \beta)$ is unitary if and only if $E_\star(\alpha, \beta)$ is a \emph{stable} parabolic bundle of parabolic degree zero. In this case, the parabolic bundle $E_\star(\alpha, \beta)$ corresponds to $\mathcal{H}_n(\alpha, \beta)$ under the Mehta-Seshadri correspondence. If $\sum(\beta_i - \alpha_i) < 0$, then $\mathcal{H}_n(\alpha, \beta)$ is unitary if and only if $E_\star(1 - \alpha, 1 - \beta)$ is a \emph{stable} parabolic bundle of parabolic degree zero. In this case, the local system $\mathcal{H}_n(\alpha, \beta)$ corresponds to the parabolic dual $E_\star(1 - \alpha, 1 - \beta)^\vee$ (equivalently, $\mathcal{H}_n(1 - \alpha, 1 - \beta)$ corresponds to $E_\star(1 - \alpha, 1 - \beta)$).
\begin{proof}
In either case, if $E_\star(\alpha, \beta)$ or $E_\star(1 - \alpha, 1 - \beta)$ is a stable parabolic bundle, then it corresponds to some unitary local system $\bl$ under the Mehta-Seshadri correspondence with the same local monodromy data at the punctures as $\mathcal{H}_n$. Since $\mathcal{H}_n$ is rigid, the local system $\bl$ must be isomorphic to $\mathcal{H}_n(\alpha, \beta)$ or $\mathcal{H}_n(1 - \alpha, 1 - \beta)$. 

Conversely, if $\mathcal{H}_n(\alpha, \beta)$ is unitary, then so is its dual $\mathcal{H}_n(1 - \alpha, 1 - \beta)$. Therefore, there are stable parabolic bundles $F_\star$ and $G_\star$ corresponding to $\mathcal{H}_n(\alpha, \beta)$ and $\mathcal{H}_n(1 - \alpha, 1 - \beta)$, respectively. If $\sum(\beta_i - \alpha_i) \geq 0$, then $\sum(\beta_i - \alpha_i) < 1$ because otherwise the sum of the parabolic weights on $F_\star$ will be less than $n$ implying that $F_\star$ is not parabolically stable. Since the sum of the parabolic weights is $n$, stability of $F_\star$ forces the underlying bundle to be isomorphic to $\mathscr{O}(-1)^{\oplus n}$. Stability is an open condition, so if $F_\star$ is stable then so is $E_\star$ as the parabolic flags of $E_\star$ are chosen generically. Rigidity of $\mathcal{H}_n$ implies that $F_\star \cong E_\star$. If $\sum(\beta_i - \alpha_i) < 0$ then $G_\star \cong E_\star(1 - \alpha, 1 - \beta)$ by the same argument. Finally, if $\sum(\beta_i - \alpha_i) = 0$ then the local monodromy at $1$ is not diagonalizable. Therefore, $\mathcal{H}_n(\alpha, \beta)$ is not unitary.
\end{proof}
\end{theorem}
Without loss of generality, we assume that the quantity $\sum (\beta_i - \alpha_i)$ is nonnegative. If it were negative, we could replace $\mathcal{H}_n(\alpha, \beta)$ in this section with its dual local system $\mathcal{H}_n(1 - \alpha, 1 - \beta)$.

\begin{definition}\label{def:interlacing}
Let $\{a_1,\dots, a_n\}$ and $\{b_1,\dots, b_n\}$ be sets of numbers on the unit circle. Suppose that $a_j = \exp(2\pi i \alpha_j)$ and $b_j = \exp(2\pi i \beta_j)$. Suppose that $0 \leq \alpha_1 < \dots < \alpha_n < 1$ and $0 \leq \beta_1 < \dots < \beta_n < 1$. We say that the sets $\{a_1,\dots, a_n\}$ and $\{b_1,\dots, b_n\}$ interlace if and only if either
\begin{enumerate}
	\item $\alpha_1 < \beta_1 < \alpha_2 < \beta_2 < \dots < \alpha_n < \beta_n$,
	\item or $\beta_1 < \alpha_1 < \beta_2 < \alpha_2 < \dots < \beta_n < \alpha_n$.
\end{enumerate}
\end{definition}

\begin{theorem}\label{thm:unitary iff interlacing}
Let $A = \{e^{2\pi i \alpha_1},\dots, e^{2\pi i \alpha_n}\}$ and $B = \{e^{2\pi i \beta_1},\dots, e^{2\pi i \beta_n}\}$. If $\sum(\beta_i - \alpha_i) \geq 0$, then $E_\star(\alpha, \beta)$ is stable and of parabolic degree zero if and only if $A$ and $B$ interlace. If $\sum(\beta_i - \alpha_i) < 0 $, then $E_\star(1 - \alpha, 1 - \beta)$ is stable and of parabolic degree zero if and only if $A$ and $B$ interlace.
\begin{proof}
Let $V = H^0 (\shom (\mathscr{O}(-1), E))$. For all $p \in \bp^1$, the fiber $E_p$ can be canonically identified with $V$ as maps $\mathscr{O}(-1) \to E$ are determined by their behavior on a single fiber. By passing through this canonical isomorphism, we obtain three distinct generic flags $T^\bullet$, $S^\bullet$, and $L^\bullet$ on $V$ from the flags on $E_0$, $E_1$, and $E_\infty$. Furthermore, a subspace $W \hookrightarrow V$ of dimension $k$ determines a map $\mathscr{O}(-1)^{\oplus k} \hookrightarrow E$ given by $W \otimes_\bc \mathscr{O}(-1) \cong \mathscr{O}(-1)^{\oplus k} \hookrightarrow V \otimes_\bc \mathscr{O}(-1) \cong E$. We can check which parabolic weights the subbundle $\mathscr{O}(-1)^{\oplus k} \hookrightarrow E$ attains by checking how the subspace $W$ intersects the various pieces of the flags $T^\bullet$, $S^\bullet$, and $L^\bullet$. Note that $\codim (T^j) = \codim(S^j) = j$ and $\codim(L^1) = n - 1$. Since the flags are generic, $\codim (T^j \cap S^k) = j + k$ and $\codim (S^j \cap L^1) = \codim(T^j \cap L^1) = n$ unless $j = 0$.

We assume $\sum(\beta_i - \alpha_i) \geq 0$ and we show $E_\star = E_\star(\alpha, \beta)$ is stable if and only if condition (1) in \autoref{def:interlacing} is satisfied. In the other case where $\sum(\beta_i - \alpha_i) < 0$, stability of $E_\star(1 - \alpha, 1- \beta)$ is equivalent to condition (2) in \autoref{def:interlacing}. The same proof applied to $E_\star(1 - \alpha, 1- \beta)$ suffices.

\textbf{Stability implies interlacing:}
To show the desired interlacing, we show two types of inequalities. First, we show that $\alpha_j < \beta_j$ for all $j$. Then we show $\beta_j < \alpha_{j + 1}$ for all $j$.

Let $F^k = \mathscr{O}(-1)^{\oplus k}$. We let $f: F^1 \hookrightarrow E$ be the map so that at $1$, the image of $f$ does not intersect $L^1$; at $0$ and $\infty$ the image intersects $T^{j - 1}$ and $S^{n - j}$, respectively. This map exists because $\codim(T^{j - 1} \cap S^{n - j}) = n - 1$ and any nonzero element in the line $T^{j-1} \cap S^{n - j} \subseteq V$ determines a map. Then, stability of $E_\star$ implies that as a parabolic subbundle of $E_\star$ the bundle $F^1_\star$ has parabolic degree smaller than zero. Therefore,
\begin{align*}
	\pardeg F^1_\star = \deg F^1 + \alpha_j + (1 - \beta_j) = \alpha_j - \beta_j < 0.
\end{align*}
	 
We now show that $\beta_j < \alpha_{j + 1}$. Let $j > 0$. We note that a map $g: F^{n - 1} \hookrightarrow E$ is determined by the choice of a subspace $W \subseteq V$ of dimension $n - 1$. Set 
\begin{align*}
W_{j + 1} := T^{j + 1} \oplus S^{n - j + 1} \oplus L^1.
\end{align*}
We let $g$ be the map defined by $W_{j + 1}$. Then, $F^{n-1}_\star$ attains the weight $\sum(\beta_i - \alpha_i)$ at $1$, the weights $\alpha_1, \alpha_2,\dots, \widehat{\alpha_{j + 1}}, \dots, \alpha_n$ at $0$, and the weights $1 - \beta_n, 1 - \beta_{n - 1}, \dots, \widehat{1 - \beta_j}, \dots, 1 - \beta_1$ at $\infty$. Here, the $\widehat{\star}$ symbol means we omit $\star$. 

Computing the parabolic degree gives us
\begin{align*}
 	\pardeg F^{n - 1}_\star = \deg F^{n - 1} + \sum_{k \neq j + 1}\alpha_k + \sum_{k \neq j} (1 - \beta_k) + \sum (\beta_i - \alpha_i)=  \beta_j - \alpha_{j + 1} < 0.
 \end{align*}

\textbf{Interlacing implies stability:}
Let $G^k_\star \hookrightarrow E_\star$ be an arbitrary parabolic subbundle of rank $k$ where $0 < k < n$. Suppose the parameters $\alpha$ and $\beta$ interlace, so that $\alpha_1 < \beta_1 < \alpha_2 < \beta_2 < \dots < \alpha_n < \beta_n$. Then, the parabolic weight attached to $L^1$ is given by $\sum (\beta_i - \alpha_i)$ and $0 < \sum (\beta_i -  \alpha_i) < 1$. We show that $\pardeg G^k_\star < 0$. We prove this in the case where $\dim(G^k_1 \cap L^1) = 0$. The other case is completely analogous.

Observe that at $0$, $G^k_\star$ attains a subset of the $\alpha_j$ as parabolic weights, and the total contribution to $\pardeg G^k_\star$ from the parabolic weights at $0$ is given by
\begin{align*}
k \alpha_1 + \sum_{j = 1}^{n - 1} (\alpha_{j + 1} - \alpha_j) \dim(G^k_0 \cap T^j).
\end{align*}
Similarly, the total contribution to $\pardeg G^k_\star$ from the parabolic weights at $\infty$ is given by
\begin{align*}
k(1 - \beta_n) + \sum_{j = 1}^{n - 1} (\beta_{j + 1} - \beta_j) \dim(G^k_\infty \cap S^{n - j}).
\end{align*}
Therefore,
\begin{align*}
\pardeg G^k_\star &= \deg G^k + k \alpha_1 + \sum_{j = 1}^{n - 1} (\alpha_{j + 1} - \alpha_j) \dim(G^k_0 \cap T^j) + k(1 - \beta_n) + \sum_{j = 1}^{n - 1} (\beta_{j + 1} - \beta_j) \dim(G^k_\infty \cap S^{n - j})\\
&= (\deg G^k + k )(1 + \alpha_1 - \beta_n) + \deg G^k(\beta_n - \alpha_1)  + \sum_{j = 1}^{n - 1} (\alpha_{j + 1} - \beta_j + \beta_j - \alpha_j) \dim(G^k_0 \cap T^j)\\
& + \sum_{j = 1}^{n - 1} (\beta_{j + 1} -\alpha_{j + 1} + \alpha_{j + 1} - \beta_j) \dim(G^k_\infty \cap S^{n - j})\\
&=(\deg G^k + k)(1 + \alpha_1 - \beta_n) + \sum_{j = 1}^{n - 1} (\alpha_{j + 1} - \beta_j) (\deg G^k + \dim (G^k_0 \cap T^j) + \dim (G^k_\infty \cap S^{n - j}))\\
& + \sum_{j = 1}^{n} (\beta_j - \alpha_j)(\deg G^k + \dim (G^k_0 \cap T^j) + \dim (G^k_\infty \cap S^{n - j + 1})).
\end{align*}

By \autoref{lem:inequalities}, we have that
\begin{align*}
\deg G^k + \dim (G^k_0 \cap T^j) + \dim(G^k_\infty \cap S^{n - j}) &\leq 0,\\
\deg G^k + \dim (G^k_0 \cap T^j) + \dim(G^k_\infty \cap S^{n - j + 1}) + \dim (G^k_1 \cap L^1) &\leq 0.
\end{align*}
Furthermore, $\deg G^k \leq -k$ and the parameters interlace so $\alpha_{j + 1} > \beta_j$, $\beta_j > \alpha_j$, and $1 + \alpha_1 - \beta_n > 0$.  Therefore,
\begin{align*}
&(\deg G^k + k)(1 + \alpha_1 - \beta_n) + \sum_{j = 1}^{n - 1} (\alpha_{j + 1} - \beta_j) (\deg G^k + \dim (G^k_0 \cap T^j) + \dim (G^k_\infty \cap S^{n - j}))\\
& + \sum_{j = 1}^{n} (\beta_j - \alpha_j)(\deg G^k + \dim (G^k_0 \cap T^j) + \dim (G^k_\infty \cap S^{n - j + 1})) < 0.
\end{align*}
This proves $\pardeg G^k_\star < 0$.
\end{proof}
\end{theorem}

We now use \autoref{thm:unitary iff interlacing} to prove the theorem of Beukers and Heckman \cite[Theorem 4.8]{Beukers_Heckman} determining when $\mathcal{H}_n$ has finite monodromy. The parameters must be rational because the local monodromy at $0$ and $\infty$ must have finite order so their eigenvalues $e^{2\pi i \alpha_j}$ and $e^{-2\pi i \beta_j}$ must be roots of unity.

\begin{corollary}\label{cor:finite iff interlacing}
The local system $\mathcal{H}_n(\alpha, \beta)$ has finite monodromy if and only if all of the $\alpha_i$ and $\beta_i$ are rational numbers and for every $\sigma \in \gal(\bq(e^{2\pi i \alpha_1}, \dots, e^{2\pi i \alpha_n}, e^{2\pi i \beta_1}, \dots, e^{2\pi i \beta_n})/\bq)$, the sets $\{\sigma(e^{2\pi i \alpha_1}), \dots , \sigma(e^{2\pi i \alpha_n}) \}$ and $\{ \sigma(e^{2\pi i \beta_1}),\dots, \sigma(e^{2\pi i \beta_n}) \}$ interlace.
\begin{proof}
	If $\mathcal{H}_n$ has finite monodromy, then every Galois conjugate of $\mathcal{H}_n$ has finite monodromy as well. Therefore, applying \autoref{thm:unitary iff interlacing} to the Galois conjugates of $\mathcal{H}_n$ gives us the desired interlacing.
	
	Conversely, suppose the interlacing condition is satisfied for all Galois conjugates. For an element $\sigma$ in the Galois group, let $\sigma \mathcal{H}_n$ be the corresponding Galois conjugate of $\mathcal{H}_n$. Applying \autoref{thm:unitary iff interlacing} to all of the Galois conjugates $\sigma \mathcal{H}_n$ implies that they are all unitary. Let $K = \bq(e^{2\pi i \alpha_1}, \dots, e^{2\pi i \alpha_n}, e^{2\pi i \beta_1}, \dots, e^{2\pi i \beta_n})$. We define the local system $\bv$ to be
	\begin{align*}
		\bv &= \prod_{\sigma \in \gal(K/\bq)} \sigma \mathcal{H}_n.
	\end{align*}
	Since $\mathcal{H}_n$ is defined over the ring of integers $\mathcal{O}_K$ (\autoref{prop:defined over ring of integers}) and the image of $\mathcal{O}_K$ is a discrete lattice inside $\prod_{\sigma \in \gal(K/\bq)} \bc$ under the diagonal embedding, the monodromy group of $\bv$ is contained in $\gl_{n \cdot [K:\bq]}(\bz)$. Furthermore, $\bv$ is unitary, so its monodromy group is a discrete subgroup of a compact group. This implies $\mathbb{V}$ has finite monodromy. Since $\mathcal{H}_n$ is a summand of $\bv$, this gives finite monodromy for $\mathcal{H}_n$ as well.
\end{proof}
\end{corollary}

\bibliographystyle{alpha}
\bibliography{hypergeometric}

%
%
%
%
%
%
%
%

\end{document}